\documentclass[12pt,letterpaper, reqno]{amsart}
\usepackage{amssymb,latexsym, amsmath, amsxtra,qtree}
\usepackage[dvips]{graphics}

\newtheorem{remark}{Remark}
\theoremstyle{plain}

\theoremstyle{definition}

\theoremstyle{remark}

\begin{document}

\title{Smooth neighbors}
\author {J. B. Conrey, M. A. Holmstrom, T. L. McLaughlin}
\begin{address}{American Institute of Mathematics\\
360 Portage Ave\\
Palo Alto, CA 94306}
\end{address}
\maketitle
\begin{abstract}
  We give a new algorithm that quickly finds smooth neighbors.
\end{abstract}
\section{Introduction}
We say that a  number is $z$-smooth if none of its prime factors exceed $z$. 
In this paper we search for solutions $b$ of 
\begin{equation} \label{eqn:main}
p\mid b(b+1)\implies p\le z.
\end{equation}
If $b$ is a solution of (\ref{eqn:main}), then we refer to the pair  $(b,b+1)$ as  $z$-smooth {\it neighbors}.
It has been known since the work of St\"{o}rmer
in 1898 [S] that for any $z$ there are only finitely many   $z$-smooth neighbors. In 1964, D. H. Lehmer [L1]
found all 869 of the 41-smooth neighbors. To do this he improved upon St\"{o}rmer's method, which relied on solving a finite number of 
Pell's equations. In fact Lehmer showed that if $b(b+1)$ is $z$-smooth then
$4b(b+1)$ is the ``$y$ part'' of the $n$th power of the fundamental solution $x_0+\sqrt{d} y_0$ of the Pell's equation
$$x^2-dy^2=1$$
where $d$ is squarefree and $z$-smooth, and where $n<(d+1)/2$. Lehmer solved all of these Pell's equations to see which led to 
41-smooth neighbors.

In 2011 Luca and Najman [LN] used a modified version Lehmer's method to find 100-smooth neighbors. In a calculation that took
15 days on a quad-core 2.66 GHz processor, they found 13,325 neighbors and claimed that this was all of the possible 100-smooth neighbors. 
The calculation was 
especially difficult because solutions to Pell's equations for squarefree 100-smooth integers can have as many as $10^{10^6}$ digits. 
Consequently they had to use a special method to even represent the solutions. In an erratum [LN1], they found 49 more solutions that 
they had missed previously.

We have found a fast, amazingly simple  algorithm that finds almost all $z$-smooth neighbors much more quickly. In fact, when we ran our method to find 100-smooth
neighbors, it completed in 20 minutes (on a similar machine as Luca and Najman used) and found 13,333 $100$-smooth neighbors.  
We were missing 37 solutions that Luca and Najman found.

Subsequently we searched for all 200-smooth neighbors. This computation took about 2 weeks and produced a list of 346,192 solutions. This 
list included all but one of the solutions from the (completed) Luca-Najman list.  We determined that this 100-smooth number missing 
from our list would be 
found using our method when searching for 227-smooth neighbors.

A plot of the logarithms of the first members from our list of  200-smooth neighbors suggests that they are normally distributed.

\section{The algorithm}

 Suppose we have a set $S$ of positive integers. For any two
 elements $b<B$ of $S$ form the ratio
 \begin{equation*}
\frac{\beta}{\beta'} =\frac{b}{b+1}\times \frac{B+1}{B}
 \end{equation*}
 where $\mbox{gcd}(\beta,\beta')=1$.
Sometimes it will be the case that $\beta'=\beta+1$, for example
$$ \frac {15}{16}=\frac 3 4 \times \frac 54.$$
We are particularly interested when this happens, i.e. we are
interested in the solutions  $\beta$ of
 \begin{equation}
\frac{\beta}{\beta+1}=\frac{b}{b+1}\times
\frac{B+1}{B}\label{eqn:beta}
 \end{equation}
 where $b$ and $B$ are in $S$.
 Given $S$, we form a new set, $S'$, which is the union of $S$ and
all of the solutions $\beta$ to (\ref{eqn:beta}). We can repeat this
process to form $(S')'=S''$ and so on. Ultimately, by St\"{o}rmer's
theorem[S], we will arrive at a set $S^{(n)}$ (meaning $n$ iterations
of priming) for which $({S^{(n)}})'=S^{(n)}$, i.e. there are no new
solutions to be added. We let
$$\delta(S)$$
denote this set $S^{(n)}$ that can no longer be enlarged by this
process. As an example, suppose that
$$S=\{1,2,3,4,5\}.$$
Then, it is easy to check that
$$S'=\{1,2,3,4,5,8,9,15,24\}.$$
and
$$S''=\{1,2,3,4,5,8,9,15,24,80\}.$$
After that we have $S'''=S''$, so
$$\delta(\{1,2,3,4,5\})=\{1,2,3,4,5,8,9,15,24,80\}.$$
Recall that Lehmer [L1] gave a complete list of 41-smooth solutions to
(\ref{eqn:main}). In particular, we can see from his list that the
above set $\{1,2,3,4,5,8,9,15,24,80\}$ is the complete list of
5-smooth solutions to (\ref{eqn:main}).  In other words, we found all
of the $z=5$ solutions to (\ref{eqn:main}) by starting with the set
$\{1,2,3,4,5\}$ and then repeatedly adding in solutions to
(\ref{eqn:beta}). This good fortune is not always the case. For
example,
$$\delta(\{1,2,3,4,5,6,7\})=\{1, 2, 3, 4, 5, 6, 7, 8, 9, 14, 15, 20, 24, 27, 35, 48, 49, 63, 80, \
125, 224, 2400\}$$ whereas from Lehmer's table we see that the
complete set of $7$-smooth solutions to (\ref{eqn:main}) includes all
of these numbers along with 4374. However, it is the case that
$$4374 \in \{1,2,3,4,5,6,7,8,9,10\}.$$
Actually, \begin{eqnarray*} \delta(\{1, 2, 5, 6, 10\})&=&\{1,2,3, 4,
5,6,7, 8, 9,10, 11, 14, 15, 20, 21, 24, 27, 32, 35,\\
&& \qquad  44, 48, 49, 54,55, 63, 80, 98, 99, 120, 125, 175, 224,\\
&& \qquad \quad 242, 384, 440, 539, 2400, 3024, 4374, 9800\}
\end{eqnarray*}
contains all of the solutions to (\ref{eqn:main}) with $y=7$. It is
easily checked that
$$\frac{4374}{4375}=\left(\frac12\right)^2\left(\frac23\right)^{-2}\left(\frac56\right)^{-4}\left(\frac67\right)
$$
so that the fraction $\frac{10}{11}$, which corresponds to $b=10$,
does not appear; but in finding the above decomposition by tracking
the iterative procedure in computing $\delta(\{1,2,5,6,10\})$ this
fraction made an appearance and then is later canceled out.

 Looking  at
Lehmer's tables we find for all primes $p$ up to 41, with the
exception of $p=7$ and $p=41$ that
$$ z_p:=\delta(\{1,2,\dots,p\}) $$
coincides exactly with the complete set of $p$-smooth solutions of
(\ref{eqn:main}). In the case $p=41$, we found that $z_{41}$ has 890 of Lehmer's solutions; 
it is missing the largest solution, 63927525375. Note that
$$ \frac{63927525375}{63927525376}=\frac{3^3\times 5^3 \times 7^7
\times 23}{2^{13}\times 11^4\times 13\times 41}.$$ The least $n$ for
which
$$63927525375\in \delta(\{1,2,\dots ,n\})$$
is $n=52$.

We calculated $z_{199}$ in a week on  a PC using Mathematica.
It has 346192 elements. A histogram of the
logarithms of these numbers is shown in Figure 1.

\begin{figure}[htp]
\begin{center}
\scalebox{0.8}[0.8]{\includegraphics{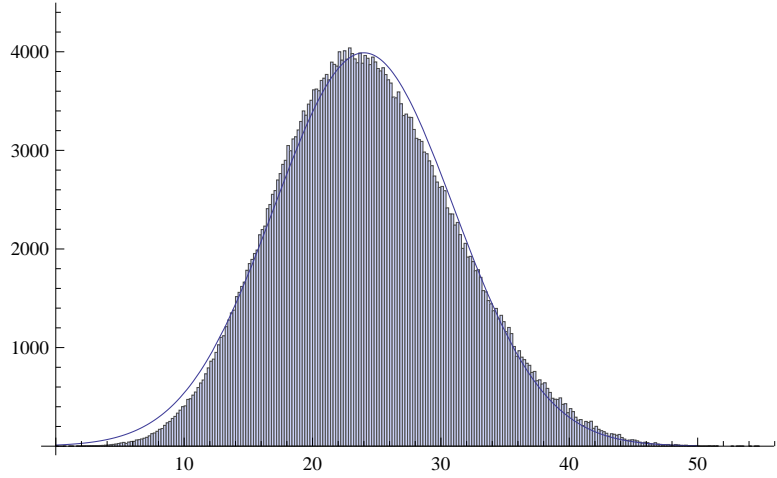}} \caption{\sf A
histogram of the set $\{\log b: p\mid b(b+1)\implies p\le 199\} $
together with the PDF of the normal distribution. }
\label{fig:lemma}
\end{center}
\end{figure}

\begin{figure}[htp]
\begin{center}
\scalebox{0.8}[0.8]{\includegraphics{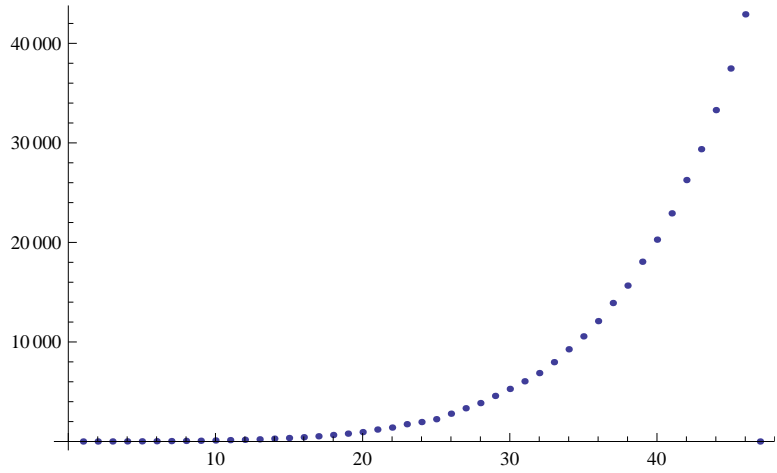}} \caption{\sf A
plot of $(n,|z_{p_n}|-|z_{p_{n-1}}|)$ where $p_n$ denotes the $n$th prime.}
\label{fig:lemma}
\end{center}
\end{figure}

 The prime before 199 is 197. The number of new $b$
in $z_{199}$, i.e. those not in
 $z_{197}$, is 43215. For all but 300 of these $b$
it is the case that $197\mid b(b+1)$. The other 300 have a smaller
largest prime factor. Here is a table of the number of these $b$ sorted by
the largest prime factor of $b(b+1)$:

\begin{eqnarray*}
\begin{array}{c|c|c|c|c|c|c|c|c|c|c|c|c|c|c|c}
\mbox{Prime}&127&131&137&139&149&151&157&163&167&173&179&181&191&193&197\\
\mbox{Number}&1&1&2&2&5&5&10&8&13&16&25&36&43&51&82
\end{array}
\end{eqnarray*}

 \section{ Solving the Diophantine equation}

In this section we analyze the solutions of
$$\frac{b}{b+1}\times \frac{B+1}{B}=\frac {\beta}{\beta+1}.$$
This equation can be written as
$$ b(B+1)(\beta+1)= (b+1)B \beta.$$
Let $\mbox{gcd}(b,B)=g$ and put
$$ b=gu\qquad \qquad B=gv$$
where clearly $\mbox{gcd}(u,v)=1$.  Then we have
$$ u (B+1) (\beta+1)= (b+1) v \beta.$$
Now, $\mbox{gcd}(u,(b+1)v)=1$ so it must be the case that $u\mid
\beta$ , say $\beta=h u$. Then we have
$$ (B+1)(\beta+1) = (b+1) v h.$$
Then $v\mid (\beta+1)$, say $\beta+1 =vx$  and $h\mid (B+1)$, say
$B+1= hy$. Thus, in terms of the variables $g,h,u,v,x,y,$ we have the
three equations
\begin{eqnarray*}
\left\{ \begin{array}{ll} hy-gv&=1\\
vx-hu&=1\\
xy-gu&=1
\end{array}\right. \end{eqnarray*}
along with $b=gu$, $B= gv$ and $\beta=hu$. If we solve for $x$ and
$y$ in the first two equations and substitute the results into the
third equation we have
\begin{eqnarray} \label{eqn:linear}
uh+gv=hv-1.\end{eqnarray} Notice that if $h$ and $v$ are given
positive coprime integers then there  exist unique positive integers
$u$ and $g$ that satisfy (\ref{eqn:linear}).  Then  $x$ and $y$ can
be determined as well as $b,B$ and $\beta$. Thus, a coprime pair $(h,v)$
leads to a triple $(b,B,\beta)$ and vice-versa.

From another perspective, one may ask how to find a pair $(b,B)$
from a given $\beta$. Take any divisor $u$ of $\beta$ and any
divisor $v$ of $\beta+1$. Then $h=\beta/u$ and $x=(\beta+1)/v$, from
which we see that $$y=\frac{v-u}{vx-hu} \qquad
g=\frac{h-x}{vx-hu}.$$
 Thus, the pair $(u,v)$ leads to the pair
$(b,B)$ given by
$$ b=\frac{uh-ux}{vx-hu}=\beta-\frac{u}{v}(\beta+1) \qquad B=\frac{hv-xv}{vx-hu}=\frac{v}{u}\beta-(\beta+1)=\frac{vb}{u},$$
where we assume that $u<v$.
And we've already seen that given $b<B$ we can define
$$ u=\frac{b}{(b,B)} \qquad v=\frac{B}{(b,B)}.$$
Thus, given $\beta$, there is a one-to-one correspondence between pairs $(u,v)$ with $u<v$, $u\mid \beta$, $v\mid (\beta+1)$ and pairs
$(b,B)$ with $b<B$ and $\frac{b}{(b+1}\times \frac{B+1}{B} =\frac{\beta}{\beta+1}.$

\section{Trees}
In this section we illustrate the way in which 601425 enters into
$z_{20}=\delta(\{1,2,\dots 20\})$.

\vspace{1.1cm}

 \Tree [.601425 [.5775 [.75 15 19 ] [.76 [.32 8 11 ]
[.56 14 19 ] ] ] [.5831 [.323 17 18 ] [.342 [.48 6 7 ]  [.56 14 19 ]
] ] ]

\vspace{1.1cm}

This tree illustrates the factorization
\begin{eqnarray*}\frac{601425}{601426}&=&\left(\frac {5775}{5776}\right)\left( \frac
{5831}{5832}\right)^{-1}\\
&=&
\left(\left(\frac{75}{76}\right)\left(\frac{76}{77}\right)^{-1}\right)\left(\left(\frac{323}{324}\right)\left(\frac
{342}{343}\right)^{-1}\right)^{-1}\\
&=& \dots \\
&=& \left(\frac  67\right)\left(\frac 7 8 \right)^{-1}\left(\frac 8
9 \right)^{-1}\left(\frac {11}{12}\right)\left(\frac
{15}{16}\right)\left(\frac {17}{18}\right)^{-1}\left(\frac
{18}{19}\right)\left(\frac {19}{20}\right)^{-1}.
\end{eqnarray*}

\section{Luca and Najman's work}

Luca and Najman published a paper in {\it Mathematics of
Computation} in 2011 which improved on Lehmer's previous work [LN]. They
gave a  list of 13325 pairs of numbers which were 97-smooth
neighbors and claimed that this list was complete.
 Their computation took 15 days on a computer with a 2.5 GHz processor.   By comparison,
when we searched for 97-smooth neighbors, our method took just under
20 minutes and, much to our surprise, produced 13333 neighbors. In
comparing their list with ours, Luca and Najman found 37 solutions
not on our list. However, we found 45 solutions not on their list!
Independently of our work, Luca and Najman found the 49 solutions they had missed previously.

Our calculations of 199-smooth numbers took more than a week and
produced 345192 solutions, some of which were the solutions we had
missed in our 97-smooth calculation.
  In fact, our $z_{199}$ is only missing one solution from the Luca-Najman list, 9591468737351909375.

\section{ The largest Luca-Najman solution}

The 79-smooth neighbor pair which starts with
$$\beta=9591468737351909375$$
was found by Luca and Najman but is not on our list (which we call
$z_{199}$) of 199-smooth neighbors. We can prove that it will first
appear in our method when we search for 227-smooth neighbors.

Given a number $\beta$, there is a 1-1 correspondence between pairs
$(b,B)$ with $b(B+1)/((b+1)B)=\beta/(\beta+1)$ and pairs $(u,v)$ with
$u\mid \beta$ and $v\mid (\beta+1)$. In one direction this
correspondence is given by
\begin{eqnarray*}
b&=& \beta-\frac{u}{v}(\beta+1)\\
 B&=& \frac{v}{u}\beta-(\beta+1)=\frac{vb}{u}
 \end{eqnarray*}

There are 1440 divisors of $\beta$ and 5632 divisors of $\beta+1$.
This means that there are a total of 8110080 pairs $(u,v)$ to
consider. For each of these pairs we computed the pair $(b,B)$ and
then computed the largest prime factor of $b(b+1)B(B+1)$.  The least
out of all of these largest prime factors was $p=227$ which appeared
for several pairs $(b,B)$, in particular for the pair
$$b = 285406166331883519 \qquad  B= 294159243066390624.$$
Therefore, $\beta$ cannot appear as a solution in our method
 unless we go up to $p=227$.

 Further analysis led us to find an explicit tree for $\beta$ all of whose bottom nodes are either in $z_{199}$
or else appear within the first two iterations (which are quick to
compute) arising in the computation of $z_{227}$; thus we can show
that $\beta\in z_{227}$ without computing all of $z_{227}$.

We now describe the tree in detail. Given a $\beta$ we have described a 1-1 correspondence between pairs $(u,v)$ and pairs $(b,B)$.
For a pair $(b,B)$ let us use the notation
$$F(b,B)=\beta$$
to denote that
$$\frac{b}{(b+1)}\times \frac{(B+1)}{B}=\frac{\beta}{(\beta+1)}.$$
Thus,
$$F(285406166331883519, 294159243066390624)=9591468737351909375.$$
is the first step of our computation.
Let's use
$$b=285406166331883519;\qquad B= 294159243066390624; \qquad \beta=9591468737351909375$$
for the rest of this section, so that $F(b,B)=\beta$.
Now let
$$g_1 = 2229716045541599; \qquad  g_2 = 2247272709023744$$
and
$$
h_1 = 186642247267999999; \qquad h_2 = 510640590102749183.$$
Then
$$F(g_1,g_2)=b; \qquad F(h_1,h_2)=B.$$
Moreover,
$$g_1\in z_{199} \qquad h_2\in z_{199}$$
so that we may restrict our attention now to $g_2$ and $h_1$.
Next, we let
$$i_1 = 907177810312319; \qquad i_2 = 911608699868750$$
and
$$
j_1 = 1671690051584; \qquad j_2 = 1672934505788.$$
Then
$$F(i_1,i_2)= h_1 \qquad F(j_1,j_2)=g_2$$
and $j_2\in z_{199}$; we are left to decompose $i_1, i_2$, and $j_1$.
Let
$$k_1 = 341611712; \qquad k_2 = 341681535;
\qquad m_1 = 300775; \qquad m_2 = 301040.$$
These satisfy
$$F(k_1,k_2)= j_1 \qquad F(m_1,m_2)=k_2$$
and $k_1,m_2\in z_{199}$.
Let
$$n_1 = 3405;\qquad  n_2 = 3444; \qquad o_1 = 454;\qquad  o_2 = 524;$$
then
$$F(n_1,n_2)=m_1; \qquad F(o_1,o_2)=n_1$$
and $n_2,o_2\in z_{199}$.
Let
$$p_1 = 199802399641; \qquad p_2 = 199846415040; \qquad q_1 = 16503580; \qquad q_2 = 16504943;$$
then $F(p_1,p_2)=i_1$; $F(q_1,q_2)=p_2$; and $p_1\in z_{199}$.
Let
$$r_1 = 561824; \qquad r_2 = 581624;\qquad s_1 = 49664; \qquad s_2 = 54480;$$
then $F(r_1,r_2)= q_1$; $F(s_1,s_2)= r_1$; and $r_2,s_1\in z_{199}$.
Let
$$t_1 = 6810; \qquad t_2 = 7783; \qquad u_1 = 108033083250000; \qquad u_2 = 122557101693480;$$
then $F(t_1,t_2)=s_2$; $F(u_1,u_2)=i_2$; and $t_2,u_1\in z_{199}$.
Let
$$v_1 = 10638314820; \qquad v_2 = 10639238337; \qquad w_1 = 1451240; \qquad w_2 = 1451438$$
then $F(v_1,v_2)=u_2$; $F(w_1,w_2)=v_1$; and $v_2,w_1\in z_{199}$.
Let
$$x_1 = 92852; \qquad x_2 = 99198; \qquad y_1 = 4312; \qquad y_2 = 4508;$$
then $F(x_1,x_2)= w_2$; $F(y_1,y_2)=x_2$; and $x_1,y_2\in z_{199}$.
Finally, $o_1, t_1, y_1 \in z_{227}$. This last fact is determined by computing the first two iterations
in the process for determining $z_{227}$ which only takes a few seconds.
\section{Computational time}

We have computed 346192 solutions to
$$ p\mid b(b+1) \implies p\le 199.$$
199 is the 46th prime. Would it be possible to use the Pell's equation method to find the complete list of 199-smooth neighbors?
No.  To use the Pell's equation method [L1] to find all of the solutions, one would
first have to find the fundamental solutions of
$$2^{46}-1 = 70368744177663 $$
different Pell's equations and then check up to the 100th solution of
each. The difficulty in finding the continued fraction expansion of
$\sqrt{d}$ grows dramatically with $d$.
We expect that for the
prime 199 we would find $d$ such that the period of $\sqrt{d}$ is as
large as $10^{40}$. Even using the subexponential algorithm employed by Luca and Najman [LN]
it would still be impossible to go this far.

Would it be possible by this method to find all of the solutions we found? No. 
For the numbers $b$ on our list $z_{199}$ the continued fraction expansions of $\sqrt{b(b+1)}$  have 
unusually small periods; in fact the largest such period is only 38.
Suppose one does a calculation where they check all $2^{46}$ square-free $d$'s that are 200-smooth to see which lead
to short continued fraction expansions as described above; 
even this calculation we estimate would take, at a conservative minimum, at least 3000 years on the computer we used.

\section{Largest solutions}
Here is a list of $(q,b)$ where $q$ is a prime number up to 197 and $b$ is the largest element of $z_{199}$
for which $q\mid(b(b+1)$ and  $p\mid b(b+1) \implies p\le q$.

\begin{eqnarray*}
&&
(2, 1), (3, 8), (5, 80), (7, 4374), (11, 9800), (13, 123200), (17, 336140),(19, 11859210),\\
&&
(23, 5142500), (29, 177182720), (31, 1611308699),(37, 3463199999), (41, 63927525375), \\
&&
(43, 421138799639), (47, 1109496723125),(53, 1453579866024), (59, 20628591204480), \\
&&
(61, 31887350832896),(67, 12820120234375), (71, 119089041053696),\\
&&  (73, 2286831727304144),(79, 1383713998733898), (83, 17451620110781856),\\
&&
 (89, 166055401586083680),(97, 49956990469100000), (101, 4108258965739505499),\\
&&
 (103, 19316158377073923834000),(107, 386539843111191224),\\
&&
 (109, 90550606380841216610), (113, 205142063213188103639),\\
&&
(127, 20978372743774437375), (131, 1043073004436787852800),\\
&&
 (137, 65244360004072055000),(139, 152295745769656587384),\\
&&
 (149, 6025407960052311234299), (151, 1801131756071318295624),\\
&&
(157, 277765695034772262487), (163, 1149394259345749379424),\\
&&
 (167, 2201197005772848768608),(173, 4574658033790609920000), \\
&&
(179, 9021820053747825025975), (181, 13989960217958128903124),\\
&&
(191, 75121996591287627735039), (193, 444171063468653314858175),\\
&&
 (197, 25450316056074220028640),(199, 589 864 439 608 716 991 201 560)
\end{eqnarray*}
Note that the Luca-Najman number 9591468737351909375 would correspond to $q=79$ above, where we have 1383713998733898.

\section{Smoothness and the ABC-conjecture}
The ABC
equation is
\begin{equation}A+B=C. \label{eqn:abc}
\end{equation}
The 1985 ABC-conjecture of Masser and Oesterl\'{e}  asserts that for any $\epsilon>0$
there is a $\kappa(\epsilon)>0$ such that for all solutions to
(\ref{eqn:abc}) the inequality
$$C < \kappa(\epsilon) \mbox{rad}(ABC)^{1+\epsilon}$$\
holds,
where $\mbox{rad}(n):=\prod_{p\mid n} p$, called the {\it radical}
of $n$, is the product of the prime divisors of $n$. In studies of
the ABC conjecture solutions to (\ref{eqn:abc}) there are a variety of indicators used to measure 
the solutions. A chief one is the 
 {\it quality}:
$$q(A,B,C)= \frac{\log C}{\log \mbox{rad}(ABC)}.$$
In fact, it seems to be standard 
to call any triple $(A,B,C)$ whose quality is at least 1 ``an ABC-triple.''
Much work has been done to systematically  find all ABC-triples. 
The first $22763667$  ABC-triples are available for
download from the abc@home web-site [A]. These include all 
ABC-triples with $C< 10^{18}$. 

Lagarias and Soundararajan [LS] define what they call the {\it smoothness exponent} $\kappa_0$ of an ABC solution by
$$\kappa_0(A,B,C) = \frac{\log z}{\log \log C}$$
where $z$ is the largest prime factor of $ABC$, and they define
 the number $\kappa_0$ by 
$$ \kappa_0:=\liminf_{C\to\infty} \kappa_0(A,B,C).$$
They observe that the ABC conjecture implies that $\kappa_0\ge 1$
and conjecture that $\kappa_0=3/2.$
 They also prove that the Generalized Riemann Hypothesis
implies that $\kappa_0\le 8$.
 
 We introduce another measure of the size  of a
$z$-smooth solution of (\ref{eqn:abc}) with a quantity we call the
{\it smoothness index}  given by
$$s(A,B,C):=\frac{\log C}{\log z}$$
where we take $z$ to be the largest prime factor of $ABC$. For
example,
$$2+25=27$$
has $C=27$ and the largest prime factor of $2\times 25\times
27=2\times 3^3 \times 5^2$ is 5 so that $s(2,25,27)=\frac{\log 27}{\log
5}=2.04782.$

{ \remark For comparison, it is easy to make the smoothness exponent large (just take $C$ to be a large prime) but hard to make it small.
For the smoothness index, the challenge is to make it large.}
{\remark The quality is more a measure of the average size of the prime
divisors of $ABC$ whereas the smoothness index measures the size of
the maximum prime divisors of $ABC$.   Note that an ABC-triple can have a high quality even if a few of the prime divisors
are large whereas that cannot happen for measures of smoothness.}
{\remark  The ABC
conjecture asserts that $q\to 1$ as $C\to \infty$. However, by work of Balog and Wooley [BW] $s(A,B,C)$ is
  unbounded.  }
 {\remark Work of Lagarias and Soundararajan [LS],
leads us to expect triples with $s(A,B,C)$ about as large $\frac{z^{2/3}}{\log z}$.
  }

\vspace{.1cm}

Thus, theoretically we know $s(A,B,C)$ can be arbitrarily large; for us the challenge is actually {\it finding} $(A,B,C)$ with $s(A,B,C)$ large.
 In other words it's a computational challenge rather than a theoretical
challenge to find large values of $s(A,B,C)$.
 Our  algorithm  finds large values of
$$s(1,B,B+1).$$
The largest value we found is
for 
$$B=19316158377073923834000$$
we have
$$s(1, B,B+1)=11.0719$$
The value of $B$ here is $\approx 1.9\times 10^{22}$ which is beyond
the range where systematic study of the ABC conjecture has been
conducted.  We know of no larger values of $s$.

\section{Some data about smoothness}
We call a triple $(A,B,C)$ {\it maximally smooth} if
$s(A_1,B_1,C_1)< s(A,B,C)$ for all $C_1< C$ (or for $C_1=C$ and
$A_1<A$). Here are the first thirteen maximally smooth triples,
together with their smoothness index:
\begin{eqnarray*}
1+3&=4 \qquad &1.262\\
3+5&=8  \qquad &1.292\\
1+8&= 9 \qquad &2\\
2+25&= 27 \qquad &2.048\\
5+27&= 32 \qquad &2.153\\
1+80&=81 \qquad & 2.730\\
3+125&=128 \qquad & 3.015\\
32+343&=375 \qquad & 3.046\\
49+576 & = 625 \qquad & 3.308\\
5+1024 &=1029 \qquad & 3.565\\
1+2400&=2401 \qquad & 4\\
1+4374 &=4375 \qquad & 4.308\\
7168+78125&=85293 \qquad & 4.427
\end{eqnarray*}
The triples on this list seem very appealing. In particular,
\begin{eqnarray*}
3+5^3&=&2^7\\
 1+5\times 2^4&=&3^4\\
 2^5+7^3 &=& 3\times 5^3\\
 5+2^{10}&=&3\times 7^3\\
 1+3\times 2^7&=&7^4\\
 7\times 2^{10} + 5^7 &=& 13\times 3^8\\
 \end{eqnarray*}
are all of the form
$$ax^\ell+by^m=cz^n$$
with $a,b,c,x,y,z$ each being 1 or a prime and $\ell,m,n>2$.

 A triple
$(A,B,C)$ is called an ``ABC-triple'' provided that $C\ge
\mbox{rad}(ABC)$. The first $22763667$  triples are available for
download from the abc@home web-site [A]. These include all of the
triples for $C< 10^{18}$. The largest $C$ for a triple on this list
is for the triple
$$ 131854322526743011 +
9091517323167918864 = 9223371645694661875$$ or in factored form
$$ 13^4\times 16651^3 + 2^4 \times 3^8 \times 53^4 \times 3313^2 =
5^4\times 7^7 \times 37 \times 59^2 \times 373^2.$$ This triple has
smoothness index $4.493$. Here we have $C= 9.2\times 10^{18}$.

Among these more than 22 million ABC - triples the triple with
largest smoothness index is
$$176202799695875+ 3178472661789594624 = 3178648864589290499$$
with smoothness index 11.0653. In factored form this is
$$\left( 5^3\times 23^3\times 42^5\right)+\left(2^{10}\times 3^7 \times 7^6 \times 13^3
\times 17^2 \times 19\right) = 11^4 \times 31 \times 37^4 \times
43^3 \times 47.$$

\section{Lehmer's table}

For each prime $z\le 41$ and each $k\le 6$ Lehmer [L2] found the
largest $n$ such that
$$p\mid n(n+1)\dots (n+h-1) \implies p\le z.$$
For example, for $z=37$ and $h=3$ he found $n=17575$ which means
that $17575\times 17576\times 17577=2^ 3\times 3^ 4\times  5^ 2
\times 7 \times  13^ 3 \times  19 \times  31 \times  37$ has no
prime factor larger than 37 and that 17575 is the largest number
with that property. Our algorithm produces many more such data
points, but in each case we can only say that our number is a {\it
lower bound} for the actual number, since we don't know for any
$z>41$ whether our list of $b$ is complete.

But to give some examples, we found that
$$ \ell(134848\times 134849 \times 134850) = 43$$
$$\ell(192510  \times 192512\times 192512 )=47$$
$$\ell(1205644 \times 1205645 \times 1205646 )= 53$$

Below is a table of lower bounds.
\begin{eqnarray*}
\begin{array}{|c|c|c|c|c|c|c|c|c|}
\hline
p& h=2&h=3&h=4&h=5&h=6&h=7\\
 \hline
 3 & 8 &  &   &   &   &     \\
 5 & 80 & 8 &   &  &  &      \\
 7 & 4374 & 48 &   &  &  &   \\
 11 & 9800 & 54 &  &  &  &   \\
 13 & 123200 & 350 & 63 & 24 &  &    \\
 17 & 336140 & 440 & 63 & 48 &  &    \\
 19 & 11859210 & 2430 & 168 & 48 &  &    \\
 23 & 11859210 & 2430 & 322 & 48 &  &    \\
 29 & 177182720 & 13310 & 322 & 54 &  &   \\
 31 & 1611308699 & 13454 & 1518 & 152 &  &    \\
 37 & 3463199999 & 17575 & 1518 & 152 &   &    \\
 41 & 63927525375 & 212380 & 1680 & 1517 & 285 &    \\
 43 & 421138799639 & 212380 & 10878 & 1517 & 340 &    \\
 47 & 1109496723125 & 212380 & 17575 & 1517 & 340 & 184   \\
 53 & 1453579866024 & 1205644 & 17575 & 1517 & 340 & 184  \\
 59 & 20628591204480 & 1205644 & 17575 & 1517 & 528 & 527  \\
 61 & 31887350832896 & 1205644 & 17575 & 1767 & 528 & 527  \\
 67 & 31887350832896 & 2018978 & 17575 & 5828 & 1271 & 527 \\
 71 & 119089041053696 & 3939648 & 17575 & 5828 & 1271 & 527  \\
 73 & 2286831727304144 & 3939648 & 70224 & 5828 & 1271 & 527  \\
 79 & 2286831727304144 & 15473808 & 70224 & 5828 & 3476 & 527  \\
 83 & 17451620110781856 & 15473808 & 97524 & 5828 & 4897 & 4896  \\
 89 & 166055401586083680 & 407498958 & 97524 & 5828 & 4897 & 4896  \\
 97 & 166055401586083680 & 407498958 & 97524 & 7565 & 7564 & 4896  \\
 \dots & \dots& \dots  &\dots  & \dots & \dots  &\dots  \\
 199& 589864439608716991201560&  768026327418&    61011223           &  1448540     &  44250    & 18904              \\
\hline
\end{array}
\end{eqnarray*}
We note some minor errors in Lehmer's tables: his $h=5$ entry for
$p=41$ should be 1517.
 Also, his Table IIA is missing the entries
10935,  12901781, and 26578125 in the 29 column and 4807, 12495,
16337, 89375 from the 31 column.

\section{Future directions}
The unanswered question in all of this is why does it work?  Can we
really get all of the $b$ for which $b(b+1)$ has only small prime
factors by repeatedly applying $\delta$ to a small initial set?

A reasonable thing to do would be to run this for a year on a
supercomputer and generate millions of smooth neighbors, perhaps
some with as many as 30 digits.

It is true that computations get longer as you raise your smoothness
barrier. But the length of computation time is due mainly to the
vast number of solutions, not the complexity in finding them.

The same technique we have elucidated above works surprisingly well
to find solutions to the more general problem
\begin{eqnarray*} \label{eqn:diff}
p\mid b(b+k) \implies p\le z \end{eqnarray*} when $k$ is an odd
integer.
 Fix a difference $k$.
Start with an initial set $S$ and let
$$ S^{'}=\{ \beta: \frac{b}{b+k}
\times \frac{B+k}{B} = \frac{\beta}{\beta+k}, \mbox{ with } b,B\in
S\};$$ then $S''=(S')'$ and so on. Iterate this procedure until it
stabilizes and call the result
$$\delta_k(S).$$
It would be interesting to do extensive computations with a range of differences $k$.

For even integers the process has to be modified somewhat. For example, to generate solutions of
$$p\mid b(b+2) \implies p\le z$$
one should start with the set $z_p=\delta_1(\{1,2,\dots,p\})$ for some $p$ and then look for solutions of
$$\frac{b}{b+1} \times \frac{B+1}{B}=\frac{\beta}{\beta+2}$$
with $b,B\in z_p$.

 \end{document}